\documentclass[a4paper]{amsart}
\usepackage{amssymb}

\usepackage{graphicx}
\usepackage[matrix,arrow,curve,tips]{xy}
           \SelectTips{cm}{}

\newtheorem{dfn}{Definition}[section]
\newtheorem{prop}[dfn]{Proposition}
\newtheorem{theo}[dfn]{Theorem}
\newtheorem{cor}[dfn]{Corollary}
\newtheorem{lem}[dfn]{Lemma}
\newtheorem{rem}[dfn]{Remark}
\newtheorem{rems}[dfn]{Remarks}
\newtheorem{ex}[dfn]{Example}
\newtheorem{exs}[dfn]{Examples}

\newcommand{\ra}{\rightarrow}
\newcommand{\lra}{\longrightarrow}
\newcommand{\oo}{\,\mbox{-}\,}
\newcommand{\com}{\mathbin{{\scriptstyle \circ }}}

\newcommand{\RR}{\mathbb{R}}

\newcommand{\cF}{\mathord{\mathcal{F}}}

\newcommand{\bY}{\bar{Y}}

\newcommand{\fb}{\mathord{\mathfrak{b}}}
\newcommand{\fk}{\mathord{\mathfrak{k}}}
\newcommand{\fg}{\mathord{\mathfrak{g}}}
\newcommand{\fh}{\mathord{\mathfrak{h}}}
\newcommand{\fl}{\mathord{\mathfrak{l}}}
\newcommand{\fX}{\mathord{\mathfrak{X}}}
\newcommand{\fL}{\mathord{\mathfrak{L}}}

\newcommand{\id}{\mathord{\mathrm{id}}}
\newcommand{\pr}{\mathord{\mathrm{pr}}}

\newcommand{\eC}{\mathord{\mathit{C}^{\infty}}}

\newcommand{\src}{\mathord{\mathrm{s}}}
\newcommand{\trg}{\mathord{\mathrm{t}}}
\newcommand{\uni}{\mathord{\mathrm{u}}}
\newcommand{\anchor}{\mathord{\mathrm{an}}}

\newcommand{\Ker}{\mathord{\mathrm{Ker}}}

\newcommand{\Aut}{\mathord{\mathrm{Aut}}}
\newcommand{\codim}{\mathop{\mathrm{codim}}}
\newcommand{\reg}{\mathord{\mathrm{reg}}}

\newcommand{\basF}{\mathord{\mathcal{F}_{\mathrm{bas}}}}

\newcommand{\baspi}{\mathord{\pi_{\mathrm{bas}}}}
\newcommand{\baspiast}{\mathord{\pi^{\ast}_{\mathrm{bas}}}}
\newcommand{\baspiinv}{\mathord{\pi^{-1}_{\mathrm{bas}}}}
\newcommand{\basOmega}{\mathord{\Omega^{0}_{\mathrm{bas}}}}

\begin{document}

\title{On the developability of subalgebroids}

\author{I. Moerdijk}
\address{Mathematical Institute, Utrecht University,
         P.O. Box 80.010, 3508 TA Utrecht, The Netherlands}
\email{moerdijk@math.uu.nl}

\author{J. Mr\v{c}un}
\address{Department of Mathematics, University of Ljubljana,
         Jadranska 19, 1000 Ljubljana, Slovenia}
\email{janez.mrcun@fmf.uni-lj.si}
\thanks{This work was supported in part by the
        Dutch Science Foundation (NWO) and
        the Slovenian Ministry of Science (M\v{S}Z\v{S} grant J1-3148)}


\begin{abstract}
In this paper, the Almeida-Molino
obstruction to developability of transversely complete
foliations is extended to Lie groupoids.
\end{abstract}

\maketitle

\section*{Introduction} \label{sec:intro}

In this paper we continue our study
--  begun in \cite{MoerdijkMrcun2004a} -- of some of the
most basic properties of Lie subgroups and subalgebras, in the wider
context of Lie groupoids and Lie algebroids. These objects occur
naturally in many contexts, such as Poisson geometry, group actions,
quantization and foliation theory
\cite{BursztynWeinstein2004,CannasdasilvaWeinstein1999,
      CattaneoFelder2000,Connes1994,
      Mackenzie1987,Pradines1967}.
One of the main features which makes
the basic
theory of Lie algebroids and Lie groupoids so much more involved
than that of Lie algebras and groups is that Lie algebroids may not be
integrable. In other words, any Lie groupoid has a Lie algebroid as its
infinitesimal part, but not every Lie algebroid arises in this way.

An obstruction to integrability
was first observed by Almeida and Molino, in the context of
foliations on compact manifolds. Recall that a foliation on a manifold
$M$ is called developable if its lift to the universal covering space
$\tilde{M}$ of $M$ is given by the fibers of a submersion
into another manifold. Given
any transversely complete foliation $\cF$ of a
compact manifold $M$, Almeida and Molino \cite{AlmeidaMolino1985}
discovered an associated Lie algebroid
$\fb(M,\cF)$, and proved that
\begin{equation} \label{eq:eq1}
\mbox{$\cF$ is developable if and only if
$\fb(M,\cF)$ is integrable.}
\end{equation}
One way to construct a developable foliation is as the kernel of a
Maurer-Cartan form on $M$ with coefficients in a Lie algebra. One can
view the Almeida-Molino result as stating that
{\em any} transversely
complete foliation on a compact manifold is the kernel of a
Maurer-Cartan form with coefficients in a Lie {\em algebroid},
and that this
Lie algebroid is integrable if and only if the foliation is developable.

Our goal is to extend this result to Lie groupoids, in the following
way: Consider a Lie groupoid $G$ over a manifold $M$,
and a subalgebroid $\fh$
of the Lie algebroid $\fg=\fL(G)$ of $G$. 
We say that $\fh$ is {\em developable} if
it can be integrated to a closed subgroupoid of the universal covering
groupoid $\tilde{G}$ of $G$. Under conditions of
source-compactness and transverse
completeness (as in the Almeida-Molino case) we will construct another
Lie algebroid $\fb(G,\fh)$, and prove an equivalence
\begin{equation} \label{eq:eq2}
\mbox{$\fh$ is developable if and only if
$\fb(G,\fh)$ is integrable.}
\end{equation}
(see Theorem \ref{theo:mcfds.25}
for a precise statement). This is formally similar to
the Almeida-Molino result (\ref{eq:eq1}),
and in fact includes the latter. To
explain the relation, recall that any manifold $M$ defines a pair
groupoid $G=M\times M$, whose universal covering groupoid
$\tilde{G}$ is the fundamental
groupoid $\Pi(M)$ given by homotopy classes of paths in $M$. Any
transversely complete foliation $\cF$ on $M$ can be viewed as a
subalgebroid of the Lie algebroid of $G$, and we will show
(Corollary \ref{cor:dbtcf.55})
that $\cF$ is a developable foliation if and only if
it can be integrated by a closed subgroupoid
of $\Pi(M)$. This shows that (\ref{eq:eq1}) is a
special case of (\ref{eq:eq2}).

In Section \ref{sec:mcfds}, we discuss Lie algebroid valued
Maurer-Cartan forms, and state our main result
(Theorem \ref{theo:mcfds.25}). For its proof,
we observe that
for a Lie groupoid $G$ over a manifold $M$,
any subalgebroid $\fh$ of the Lie
algebroid $\fL(G)$ of $G$ defines a foliation
$\cF(\fh)$ of $G$, whose leaves are
contained in the fibers of the source map
$\src\!:G\ra M$. 
In Section \ref{sec:btcf} we will
show that if the foliation $\cF(\fh)$ is
locally transversely parallelizable
and the fibers of $\src$
are compact, then $\cF(\fh)$ makes
$\src\!:G\ra M$ into a locally trivial bundle of
transversely complete
foliations. This fact is independent of the groupoid structure on $G$,
and we will prove it in the general context of
a submersion $s\!: N\ra M$
with suitably foliated fibers, see Theorem \ref{theo:btcf.11}.
Next, we will give various
characterizations of the fiberwise developability of a locally trivial
bundle $N\ra M$ of foliated manifolds, and prove that this property is
equivalent to the integrability of a specific algebroid
(Theorem \ref{theo:dbtcf.15}).
Finally, we will
apply our results on bundles of foliated manifold to groupoids,
and prove Theorem \ref{theo:mcfds.25}.

\section{Developable subgroupoids and Maurer-Cartan forms} \label{sec:mcfds}

We begin by introducing the notion of {\em developability}
for subalgebroids, which is motivated by the 
well-known notion of
developability for foliations.

For a Lie groupoid $G$ over a manifold $M$,
we shall denote by $\src\!:G\ra M$ the {\em source} map
of $G$, and by $\trg\!:G\ra M$ the {\em target} map of $G$.
Recall that $G$ is said to be
{\em source-connected},
if the fibers of the source map of $G$
are connected, and 
{\em source-simply connected} if the fibers
of the source map of $G$ are simply connected.
In this paper, we shall assume that all
Lie groupoid are source-connected.

Let $G$ be a (source-connected)
Lie groupoid over a manifold $M$.
Denote by $\fg=\fL(G)$ the Lie algebroid
associated to $G$. We denote by $\tilde{G}$ the
source-simply connected covering groupoid of $G$
which has the same Lie algebroid as $G$.

We shall consider subalgebroids of
$\fg$ over the same base $M$.
Recall from \cite{MoerdijkMrcun2004a} that
such a subalgebroid $\fh$ of $\fg$ corresponds to
a right-invariant foliation $\cF(\fh)$ of $G$
which refines the foliation $\cF(\src)$ of $G$
given by the fibers of the source map,
$\cF(\fh)\subset\cF(\src)$.
Any (injectively immersed source-connected)
Lie subgroupoid $H$ of $G$
over the same base $M$ gives rise to a Lie subalgebroid
$\fh$ of $\fg$, and in turn the leaves of $\cF(\fh)$
are precisely the right cosets of $H$ in $G$.
In this case we write $\cF(H)=\cF(\fh)$,
and we say that $H$ integrates $\fh$,
or more precisely, that the inclusion $H\ra G$ integrates
the inclusion $\fh\ra\fg$.
We denote by $G/H$ the space of leaves of the
associated foliation $\cF(H)$.

We say that  a subalgebroid $\fh$ of $\fg$ is
{\em developable}
if $\fh$ can be integrated to a closed subgroupoid of the
source-simply connected cover $\tilde{G}$ of $G$.
A subgroupoid $H$ of $G$ is {\em developable}
if the associated  Lie subalgebroid of $\fg$ is.

Recall that a foliation $\cF$ of a manifold $N$
is {\em simple} if it is given by the components of the
fibers of a submersion into a Hausdorff manifold.
It is called {\em strictly simple} if there exists a
smooth structure of a Hausdorff manifold on the space
of leaves
$N/\cF$ such that the quotient map $N\ra N/\cF$ is
a submersion.
The results of \cite{MoerdijkMrcun2004a}
provide the following equivalent formulation of
developability:

\begin{prop} \label{prop:mcfds.4}
Let $G$ be an source-connected Lie groupoid with Lie algebroid $\fg$,
and let $\fh$ be a subalgebroid of $\fg$.
The following two statements are equivalent:
\begin{enumerate}
\item [(i)]  The Lie algebroid $\fh$ is developable.
\item [(ii)] The pull-back of the foliation $\cF(\fh)$
             to the source-simply connected cover of $G$
	     is strictly simple.
\end{enumerate}
\end{prop}

\begin{proof}
This equivalence follows directly from
\cite[Proposition 3.1]{MoerdijkMrcun2004a},
\cite[Proposition 3.2 (ii)]{MoerdijkMrcun2004a}
and the fact that pull-back of the foliation $\cF(\fh)$
to the source-simply connected cover $\tilde{G}$ of $G$
is the foliation of $\tilde{G}$ given by the
same subalgebroid
$\fh\subset\fg=\fL(\tilde{G})$.
\end{proof}

\begin{ex} \rm \label{ex:mcfds.2}
Recall that a foliation $\cF$ of a manifold $N$
is called {\em developable} if
its pull-back to the universal covering space of
the manifold is strictly simple. 
Such a foliation $\cF$ can be viewed as a subalgebroid
of the Lie algebroid $T(M)$.
The Lie algebroid $T(M)$ is integrated by the
pair groupoid $M\times M$, but also by the
source-simply connected fundamental groupoid
$\Pi(M)$ of $M$. By 
Corollary \ref{cor:dbtcf.55} below,
the foliation $\cF$ is developable if
and only if $\cF$, viewed as a subalgebroid of $T(M)$,
is integrable by a closed subgroupoid of $\Pi(M)$.
In this way, our definition of developability of
subalgebroids extends the usual one for foliations.
\end{ex}

An important class of examples are the subalgebroids
which arise as kernels of Maurer-Cartan forms, which we now discuss.

Let $\fg$ be a Lie algebroid over $M$. Let $\fk$ be a finite
dimensional Lie algebra. Then $\fk$ pulls back to an algebroid
$\fk_{M}=\fk\times M$ over $M$, with zero anchor.
Suppose that we are given a (left) action $\nabla$ of $\fg$ on $\fk_{M}$
along the identity map of $M$. Recall from
\cite{HigginsMackenzie1990,MoerdijkMrcun2002}
that such an action assigns to each section $X\in\Gamma(\fg)$
a derivation $\nabla_{X}$ on the Lie algebra $\eC(M,\fk)$
of smooth $\fk$-valued functions, $\eC(M)$-linear and flat in $X$,
$$ \nabla_{[X,Y]}=\nabla_{X}\nabla_{Y}-\nabla_{Y}\nabla_{X}\;,$$
which moreover satisfies the Leibniz law
$$ \nabla_{X}(f\alpha)=f\nabla_{X}(\alpha)+\anchor(X)(f)\alpha $$
for any $f\in\eC(M)$ and any $\alpha\in\eC(M,\fk)$.

A {\em Maurer-Cartan form} on  $\fg$
with coefficients in $\fk$ is a map
$\omega\!:\fg\ra\fk_{M}$
of vector bundles over $M$
satisfying the Maurer-Cartan equation
$$ d\omega+\frac{1}{2}[\omega,\omega]=0\;.$$
Here $d$ is the differential of the Lie algebroid cohomology of $\fg$
with coefficients in $\fk$,
$$ 2d\omega(X,Y)=\nabla_{X}(\omega(Y))-\nabla_{Y}(\omega(X))-\omega([X,Y])\;,$$
while the bracket on maps
$\omega,\lambda\!:\fg\ra\fk_{M}$ of vector bundles over $M$
is given by
$$ 2[\omega,\lambda](X,Y)=[\omega(X),\lambda(Y)]
   -[\omega(Y),\lambda(X)] \;.$$

\begin{ex} \rm \label{ex:mcfds.6}
Let $\fk$ be a finite dimensional Lie algebra,
and let $\fg$ be a Lie algebroid over $M$.
The {\em trivial action} $\nabla^{\mathrm{triv}}$
of $\fg$ on $\fk\times M=\fk_{M}$,
along the identity map of $M$, is given by the derivative
along the anchor,
$$ \nabla^{\mathrm{triv}}_{X}(\alpha)=\anchor(X)(\alpha)\;.$$
For $\fg=T(M)$, a Maurer-Cartan form
on $T(M)$ with values in $\fk$ with respect to the trivial action
is a usual Maurer-Cartan form on the manifold $M$.
\end{ex}

\begin{lem} \label{lem:mcfds.7}
Let $\fg$ be a Lie algebroid over $M$, let $\fk$ be a finite
dimensional Lie algebra, and suppose that
$\fk_{M}$ is equipped with an action $\nabla$ of $\fg$
along the identity map. For a map $\omega\!:\fg\ra\fk_{M}$ of
vector bundles over $M$, the following conditions are equivalent:
\begin{enumerate}
\item [(i)]
   $\omega$ is a Maurer-Cartan form on $\fg$ with respect
   to the action $\nabla$.
\item [(ii)]
   $\omega([X,Y])=[\omega(X),\omega(Y)]
    +\nabla_{X}(\omega(Y))-\nabla_{Y}(\omega(X))$
   for any $X,Y\in\Gamma(\fg)$.
\item [(iii)]
   $(\id,\omega)\!:\fg\ra\fg\ltimes\fk_{M}$ is a morphism of
   Lie algebroids over $M$.
\end{enumerate}
\end{lem}

In (iii), the semi-direct product $\fg\ltimes\fk_{M}$ is the bundle
$\fg\oplus\fk_{M}$, with bracket defined by
$$ [(X,\alpha),(Y,\beta)]=
   ([X,Y],[\alpha,\beta]+\nabla_{X}(\beta)-\nabla_{Y}(\alpha)) $$
for any $X,Y\in\Gamma(\fg)$ and $\alpha,\beta\in\eC(M,\fk)$.
The equivalence of conditions (i)-(iii) is a trivial calculation.

\begin{prop} \label{prop:mcfds.13}
Let $\fk$ be a finite-dimensional Lie algebra,
let $G$ be a source-connected Lie groupoid
with Lie algebroid $\fg$ acting on the Lie
algebroid $\fk_{M}=\fk\times M$, and let
$\omega\!:\fg\ra\fk_{M}$ be a non-degenerated
Maurer-Cartan form on $\fg$. Then $\Ker(\omega)$
is a developable subalgebroid of $\fg$.
\end{prop}

\begin{proof}
Let $\tilde{G}$ be the source-simply connected
Lie groupoid covering $G$, and let $K$ be the
simply connected Lie group with Lie algebra $\fk$.
Write $K_{M}=K\times M$ for the trivial bundle of
Lie groups over $M$ with fiber $K$,
which integrates the Lie algebroid $\fk_{M}$.
By \cite{MoerdijkMrcun2002},
the action of $\fg$ on $\fk_{M}$ integrates to
an action of $\tilde{G}$ on $K$
(more precisely, an action on $K_{M}$), and we can
form the semi-direct product
$\tilde{G}\ltimes K_{M}$ over $M$.
(Its arrows are pairs $(\tilde{g},k)$, where
$\tilde{g}\!:x\ra y$ is an arrow in $\tilde{G}$
and $k\in K$; composition is given by
$(\tilde{h},l)(\tilde{g},k)=(\tilde{h}\tilde{g},l(\tilde{h}k))$.)
By Lemma \ref{lem:mcfds.7}
and \cite[Proposition 3.5]{MoerdijkMrcun2002},
the morphism of Lie algebroids $(\id,\omega)$
integrates to a morphism of Lie groupoids
$(\id,\Omega)\!:\tilde{G}\ra\tilde{G}\ltimes K_{M}$.
Then $\Omega\!:\tilde{G}\ra K_{M}$
is a `twisted' homomorphism,
$$ \Omega(\tilde{h}\tilde{g})=\Omega(\tilde{h})(\tilde{h}\Omega(\tilde{g}))\;,$$
whose kernel is a closed subgroupoid of $\tilde{G}$.
This closed subgroupoid integrates the Lie algebroid
$\Ker(\omega)$,
hence the latter is developable.
\end{proof}

Let $\cF$ be a foliation of a manifold $N$.
Recall that a vector field $Y$ on $N$ is {\em projectable}
with respect to $\cF$ if its (local) flow preserves the
foliation, or equivalently, if the Lie derivative
of $Y$ in the direction of a vector field tangent to
$\cF$ is again tangent to $\cF$.
A foliation $\cF$ of $N$ is {\em transversely
complete} if any tangent vector on $N$
can be extended to a complete
projectable vector field on $N$
(see also Remarks \ref{rems:btcf.1} below).

A subgroupoid $H$ of a Lie groupoid $G$ is
{\em transversely complete} if the
associated foliation $\cF(H)$ of $G$ is
transversely complete.
For example, any transitive subgroupoid $H$
of $G$ is automatically transversely complete
\cite[Proposition 3.5]{MoerdijkMrcun2004a}.

Our aim is to prove that for any transversely
complete subgroupoid $H$ of a source-compact groupoid $G$
(i.e.\ $G$ is Hausdorff and its source map is proper),
its Lie algebra $\fh$ is the kernel of a Maurer-Cartan form.
This result can be interpreted
 as a generalization of Molino's theorem
for transversely complete foliations.
As in Molino's case, it requires
the Maurer-Cartan form to take values in a Lie algebroid
(see Example \ref{exs:mcfds.20} (2)).

Let $G$ be a Lie groupoid over $M$ with Lie algebroid $\fg$.
Let $\fk$ be a Lie algebroid over $W$ equipped with an action
of $\fg$ along a submersion $f\!:W\ra M$. Recall from
\cite{MoerdijkMrcun2002} that this means in particular
that $\fk\ra W$ is a bundle of Lie algebroids over $M$,
and that $\fg$ acts on $W$ and $\fk$, via suitable maps
$R\!:\Gamma(\fg)\ra\fX(W)$ and $\nabla$ assigning to each
$X\in\Gamma(\fg)$ a derivation $(\nabla_{X},R(X))$ on $\fk$.
A {\em Maurer-Cartan form} on $\fg$ with values in $\fk$
is a section of the projection map of Lie algebroids
$\fg\ltimes\fk\ra\fg$ over $f\!:W\ra M$.
Such a section is
given by a map
$\omega\!:\fg\ra\fk$
of vector bundles over a section $\alpha\!:M\ra W$ of $f$,
with the property that
$$ ((\alpha,\id),\omega)\!:\fg\lra f^{\ast}(\fg)\oplus\fk $$
defines a morphism of Lie algebroids
$\fg\ra\fg\ltimes\fk$
over $\alpha\!:M\ra W$.

\begin{rem} \rm \label{rem:mcfds.17}
It is of course possible to spell out the last condition
in detail:
First of all,
compatibility with the anchor maps of $\fg$ and 
$\fg\ltimes\fk$ means that $\omega$ should satisfy the identity
\begin{equation} \label{eq:mcfds.18}
\anchor(\omega(v))=(d\alpha)_{x}(\anchor(v))-R_{\alpha(x)}(v)
\end{equation}
for any $x\in M$ and any $v\in\fg_{x}$;
next, compatibility of the brackets can be expressed by the
equation
\begin{equation} \label{eq:mcfds.19}
\omega\com[X,Y]=([\omega(X),\omega(Y)]
+\nabla_{X}(\omega(Y))-\nabla_{Y}(\omega(X)))\com\alpha\;.
\end{equation}
Here $X$ and $Y$ are sections of $\fg$, while
$\omega(X)$ and $\omega(Y)$ denote arbitrary
sections of $\fk$ which extend
$\omega\com X\!:M\ra\fk$ respectively $\omega\com Y\!:M\ra \fk$
along the embedding $\alpha\!:M\ra W$.
It follows from (\ref{eq:mcfds.18}) that
the right hand side of (\ref{eq:mcfds.19})
is independent of the choice of these extensions.
\end{rem}

\begin{exs} \rm \label{exs:mcfds.20}
(1)
A Maurer-Cartan form with values in a Lie algebra $\fk$
(as discussed above)
can be seen as a special case of a Maurer-Cartan form
with values in the Lie algebroid $\fk_{M}=\fk\times M$,
where the action is taken along the identity map
of the base manifold $M$.

(2) (Molino)
Let $M$ be a compact manifold equipped with a 
transversely complete foliation $\cF$.
The closures of the leaves of such a foliation
are the fibers of a locally trivial fiber bundle
$M\ra W$ of foliations \cite{Molino1977}.
There is the associated Lie algebroid
$\fb(M,\cF)$ over $W$
such that $\cF$ is the kernel
of a map of Lie algebroids over the 
projection $M\ra W$
 \cite{MoerdijkMrcun2003,Molino1988}.
\begin{equation} \label{eq:mcfds.25}
\xymatrix{
T(M) \ar[d] \ar[r] & \fb(M,\cF) \ar[d] \\
M \ar[r] & W
}
\end{equation}
To see this as a Maurer-Cartan form,
consider the Lie algebroid 
$\fb(M,\cF)\times M$ over $W\times M$.
This is a bundle of Lie algebroids over $M$,
and has the canonical action of the pair
groupoid $M\times M$
(just by acting on the second component only).
This action differentiates to an action of
$T(M)$ on $\fb(M,\cF)\times M$.
The morphism of algebroids in (\ref{eq:mcfds.25})
induces a section of the projection
$T(M)\ltimes (\fb(M,\cF)\times M)\ra T(M)$.
\end{exs}

Recall from \cite[Theorem 3.7]{MoerdijkMrcun2004a}
that if $H$ is a transversely complete subgroupoid
of a Hausdorff Lie groupoid $G$, then
its closure $\bar{H}$ in $G$ is a closed
Lie subgroupoid of $G$, and the space of right
cosets $G/\bar{H}$ has the natural structure of a
Hausdorff manifold such that the projection
$G\ra G/\bar{H}$ is a locally trivial fiber bundle.

We can now state our main theorem, which
we shall prove in Section  \ref{sec:pmt}.

\begin{theo} \label{theo:mcfds.25}
Let $G$ be an source-compact
Lie groupoid with Lie algebroid $\fg$,
and let $H$ be a transversely
complete subgroupoid of $G$ with Lie
algebroid $\fh\subset\fg$.

(i)
There exists a regular Lie algebroid
$\fb(G,\fh)$ over $G/\bar{H}$ such that
the projection $G\ra G/\bar{H}$
lifts to a natural
surjective morphism of Lie algebroids
$\cF(\src)\ra\fb(G,\fh)$ with kernel $\cF(\fh)$.

(ii)
The natural action of $\fg$ on $G/\bar{H}$ lifts
to an action of $\fg$ on $\fb(G,\fh)$.

(iii)
The Lie algebroid
$\fh$ is the kernel of a Maurer-Cartan form
on $\fg$ with values in $\fb(G,\fh)$.

(iv)
The Lie algebroid
$\fh$ is developable if and only if 
$\fb(G,\fh)$ is integrable.
\end{theo}

The assumption that $G$ is source-compact
is used to ensure that the source map
of $G$ is a locally trivial bundle of
transversely complete foliations
(Theorem \ref{theo:btcf.11}), and can
in fact be replaced by this weaker condition.

Part (ii) of Theorem \ref{theo:mcfds.25}
states in particular
that $\fb(G,\fh)$ is in fact a bundle of Lie algebroids
over the base space $M$ of $G$; more precisely,
its anchor  is annihilated by the map
$T(G/\bar{H})\ra T(M)$ induced by the differential
of the source map. Thus, $\fb(G,\fh)$
is a bundle of Lie algebras in case the map
$G/\bar{H}\ra M$, induced by the source,
is a diffeomorphism. The result of
Douady and Lazard \cite{DouadyLazard1966}
thus gives the following corollary:

\begin{cor} \label{cor:mcfds.36}
Let $G$ be an source-compact
Lie groupoid with Lie algebroid $\fg$,
and let $H$ be a transversely
complete subgroupoid of $G$ with Lie
algebroid $\fh\subset\fg$.
Then $H$ is developable whenever
the cosets of $H$ are dense in the
source-fibers of $G$.
\end{cor}

\begin{ex} \rm \label{ex:mcfds.36}
Let $\fl$ be the Lie algebra of a compact
connected Lie group $L$, and let $\omega$
be a non-singular
Maurer-Cartan form with values in $\fl$
on a connected compact manifold $M$.
Then $\omega$ defines a transitive (hence transversely complete)
subgroupoid  $H$ of the groupoid $M\times L\times M$
whose source-fibers are the holonomy covers of $\omega$,
see \cite[Example 3.8]{MoerdijkMrcun2004a}.
If the holonomy group of $\omega$ is dense in $L$,
then the cosets of $H$ are dense in $M\times L$, hence
$H$ is developable by Corollary \ref{cor:mcfds.36}.
\end{ex}

\section{Bundles of transversely complete foliations} \label{sec:btcf}

First, we recall some standard definitions.
Let $\cF$ be a foliation of a manifold $N$.
Write $\fX(\cF)$ for the Lie algebra of vector
fields on $N$ which are tangent to $\cF$, and
let $L(N,\cF)$ be the normalizer of
$\fX(\cF)$ in the Lie algebra $\fX(N)$ of vector fields
on $N$. Vector fields in $L(N,\cF)$ are called
{\em projectable} vector fields on $(N,\cF)$, and
the quotient
$$ l(N,\cF)=L(N,\cF)/\fX(\cF) $$
is the Lie algebra of transverse vector fields on $(N,\cF)$.
The transverse vector fields on
$(N,\cF)$ can be identified with the holonomy invariant
sections of the normal bundle $\nu(\cF)=T(N)/T(\cF)$.
Both $L(N,\cF)$ and $l(N,\cF)$ are modules
over the algebra $\basOmega(N,\cF)$
of basic functions on $(N,\cF)$, and the quotient projection
$$ L(N,\cF) \lra l(N,\cF) \;,\;\;\;\;\;\;\;\;
   Y\mapsto\bY\;,$$
is $\basOmega(N,\cF)$-linear.

A foliation $\cF$ of $N$ is {\em transversely parallelizable}
\cite{Conlon1974,Molino1988}
if there exists a global frame of $\nu(\cF)$ consisting
of transverse vector fields on $(N,\cF)$.
Such a frame of $\nu(\cF)$ is referred to as {\em transverse
parallelism} on $(N,\cF)$.

A foliation $\cF$ of $N$ is {\em locally transversely
parallelizable} if for any $y\in N$ and any tangent vector
$\xi\in T_{y}(N)$, there exists a projectable vector field
$Y$ on $(N,\cF)$ such that $Y_{y}=\xi$.
A foliation $\cF$ of $N$ is {\em transversely complete}
\cite{Molino1977,Molino1988}
if for any $y\in N$ and any tangent vector
$\xi\in T_{y}(N)$ there exists a complete projectable vector
field $Y$ on $(N,\cF)$ such that $Y_{y}=\xi$.

\begin{ex} \rm \label{ex:mcfds.15}
Consider a Maurer-Cartan form $\omega$ as in
Proposition \ref{prop:mcfds.13}.
Let $\fk_{G}$ be the trivial bundle of Lie algebras over $G$,
and let $\cF(\src)$ be the foliation on $G$
by the fibers of the source map.
Then $\omega$ defines a flat $\cF(\src)$-partial
connection on $\fk_{G}$, and hence a `monodromy' map
$$ \Pi^{\src}(G)\lra\Aut(\fk)\;.$$
If this map is trivial, i.e.\ if $\omega$
has trivial monodromy, then the invariant foliation
corresponding to $\Ker(\omega)$ is transversely parallelizable.
In particular, if $G$ is source-compact  and
$\omega$ has finite monodromy, then there exists
an source-compact covering groupoid $G'$ of $G$
with the same algebroid $\fg$, on which $\Ker(\omega)$
defines a transversely parallelizable foliation.
\end{ex}

\begin{rems} \rm \label{rems:btcf.1}
(1)
Note that a foliation $\cF$ of $N$ is
locally transversely parallelizable if and only if
the evaluation map $l(N,\cF)\ra \nu_{y}(N,\cF)$
is surjective for any $y\in N$.
Equivalently, a foliation $\cF$ of $N$
is locally transversely parallelizable if for any
$y\in N$ there exist projectable vector fields
$Y_{1},\ldots,Y_{q}$ such that
$((\bY_{1})_{y},\ldots,(\bY_{q})_{y})$
is a basis of $\nu_{y}(N,\cF)$.
In particular, any transversely parallelizable foliation
is locally transversely parallelizable.

Any locally transversely parallelizable foliation
of a compact manifold is transversely complete.
Examples of transversely complete foliations also
include foliations given by the fibers of locally trivial fiber bundles,
and Lie foliations on compact manifolds \cite{Fedida1971}.

(2)
Let $\cF$ be a foliation on $N$, and let
$\bY_{1},\ldots,\bY_{q}$ be transverse vector fields on
$(N,\cF)$. We say that $y\in N$ is a {\em regular point}
of $(\bY_{1},\ldots,\bY_{q})$ if the normal vectors
$(\bY_{1})_{y},\ldots,(\bY_{k})_{y}\in \nu_{y}(N,\cF)$
form a basis of $\nu_{y}(\cF)$. 
The {\em regular set} of
$(\bY_{1},\ldots,\bY_{q})$ is the set of all regular points
of $(\bY_{1},\ldots,\bY_{q})$, and will be denoted by
$$ \reg(\bY_{1},\ldots,\bY_{q})\subset N \;.$$
This set is open in $N$, and also $\cF$-saturated
because of the holonomy invariance of the transverse
vector fields.
We say that $(\bY_{1},\ldots,\bY_{q})$ is a
{\em local transverse parallelism} on $(N,\cF)$ with regular
set $\reg(\bY_{1},\ldots,\bY_{q})$.

Note that $\cF$ is locally transversely parallelizable
if and only if
the regular sets of local transverse parallelisms on $(N,\cF)$
cover $N$.
If $\cF$ is a locally transversely parallelizable foliation of
a manifold $N$, then any leaf $L$ of $\cF$
has an open $\cF$-saturated neighbourhood $U\subset N$ on which
the foliation $\cF|_{U}$ is transversely parallelizable.

(3)
Any locally transversely parallelizable foliation has
trivial holonomy. Furthermore,
any transversely complete
foliation $\cF$ of a connected manifold $M$
is homogeneous, i.e.\ the group $\Aut(M,\cF)$ of its
foliation automorphisms acts transitively on $M$. 
This is a direct consequence of the fact that
a  vector field is projectable if and only if its
(local) flow preserves
the foliation.

Any strictly simple foliation $(M,\cF)$
is simple; if the foliation is homogeneous,
then these two notions coincide
\cite[Theorem 4.3 (vi)]{MoerdijkMrcun2003}.
Note that simple
foliations are preserved under the
pull-back along a covering projection,
while strictly simple are not.
Since the homogeneity is also preserved
under the pull-back along a covering projection,
a homogeneous foliation is developable if
and only if its pull-back to a covering
space of the manifold is simple.

(4)
Recall that a homogeneous foliation $\cF$ of $N$
admits an associated {\em basic} foliation $\basF$, given by
$$ \fX(\basF)=\{ X\in\fX(N) \,|\, X(\basOmega(N,\cF))=0 \}
   \supset \fX(\cF) \;. $$
The foliation $\basF$ is again homogeneous, and satisfies
$\basOmega(N,\basF)=\basOmega(N,\cF)$
and $L(N,\cF)\subset L(N,\basF)$. Furthermore,
the space of basic leaves $N/\basF$ has a natural structure
of a Hausdorff manifold such that the basic projection
$\baspi\!:N\ra N/\cF$ is a submersion, and
$\basOmega(N,\cF)=\eC(N/\basF)$.

(5)
Suppose that $\cF$ is a transversely complete foliation
of a manifold $N$.
By Molino's structure theorem \cite{Molino1977},
the closures of the leaves of a transversely complete
foliation $\cF$ of $M$ are the fibers of
the associated basic fibration $N\ra N/\basF$,
which is in fact a locally
trivial fiber bundle of (Lie) foliations.
In particular, this implies that 
if two transversal vector fields on $(N,\cF)$
agree at a point of $N$, they also agree along
the basic leaf through that point.
Furthermore, the regular set of any
local transverse parallelism on $(N,\cF)$
is $\basF$-saturated.
\end{rems}

Let $\cF$ be a locally transversely parallelizable
foliation of a connected manifold $N$, and
let $s\!:N\ra M$ be a
surjective submersion with connected compact
$\cF$-saturated fibers.
Since the image of the homomorphism
$s^{\ast}\!:\eC(M)\ra\eC(N)$ is a subalgebra of
$\basOmega(N,\cF)$ and the fibers of $s$
are compact, it follows that $\cF$ is
transversely complete and therefore homogeneous.
It is clear that the compact fibers of $s$ are also
$\basF$-saturated.
Furthermore, the submersion $s$ factors through
the associated basic fibration
$\baspi\!:N\ra N/\basF$
as a surjective submersion
$\bar{s}\!:N/\basF\ra M$
with compact connected fibers.

The aim of this section is to
prove the following theorem, which says that
any such surjective 
submersion $s\!:(N,\cF)\ra M$
with connected compact
$\cF$-saturated fibers
is actually a locally trivial bundle of
foliated manifolds.

\begin{theo} \label{theo:btcf.11}
Let $\cF$ be a locally transversely parallelizable
foliation of codimension $q$
of a connected manifold $N$, and
let $s\!:N\ra M$ be a
surjective submersion with connected compact
$\cF$-saturated fibers.
Put $m=\dim(M)$ and $\bar{q}=\codim(\basF)$.
Then for any $y\in N$ there exist
vector fields $Y_{1},\ldots,Y_{q}\in L(N,\cF)\cap L(N,\cF(s))$
such that
\begin{enumerate}
\item [(i)] 
    $(\bY_{1},\ldots,\bY_{q})$ is a local transverse
    parallelism on $(N,\cF)$, regular on a neighbourhood of $y$,
\item [(ii)]
    $ds(Y_{i})=0$ for $i=m+1,\ldots,q$, and
\item [(iii)]
    $d\baspi(Y_{i})=0$ for $i=\bar{q}+1,\ldots,q$.
\end{enumerate}
In particular,
the restriction of $\cF$ to any fiber $N_{x}$ of $s$ over $x\in M$
is a transversely complete foliation, and
the submersion $s\!:(N,\cF)\ra M$ is a locally trivial
fiber bundle of foliated compact manifolds.
\end{theo}

As we shall show,
this theorem is a consequence
of Proposition \ref{prop:btcf.10} below,
which we formulate and prove first.

Suppose that $\phi\!:N\ra M$ is a surjective submersion,
and write $\cF(\phi)$ for the foliation of $N$
given by the connected components of the fibers of $\phi$.
Let $C$ be a subspace of $L(N,\cF(\phi))$ and
$D$ a subspace of $\fX(M)$ such that
$d\phi(C)\subset D$. Then we say that the linear map
$d\phi\!:D\ra C$ has the {\em lifting property}
if for any vector field
$X\in C$ and any $\xi\in T_{y}(N)$
such that $d\phi(\xi)=X_{\phi(y)}$ there exists
a vector field $Y\in D$ such that $Y_{y}=\xi$ and
$d\phi(Y)=X$.
Thus the lifting property of
$d\phi\!:D\ra C$ in particular implies that
$d\phi(D)=C$.

\begin{ex} \rm \label{ex:btcf.9}
Let $\phi\!:N\ra M$ be a surjective submersion,
and denote by $\cF(\phi)=\Ker(d\phi)$ the foliation of $N$
with the connected components of the fibers of $\phi$ as leaves.
Then the Lie algebra homomorphism
$d\phi\!:L(N,\cF(\phi))\ra \fX(M)$ has the
lifting property. 
In particular, the foliation
$\cF(\phi)$ is locally transversely parallelizable.

To see this,
take any vector field $X\in\fX(M)$ and any
$\xi\in T_{y}(N)$ with $d\phi(\xi)=X_{\phi(y)}$.
Note that we can choose
an open cover $(U^{\alpha})$
of $N$ and vector fields
$Y^{\alpha}\in L(U^{\alpha},\cF(\phi)|_{U^{\alpha}})$
such that $d\phi(Y^{\alpha})=X|_{U^{\alpha}}$.
Furthermore, we can assume that $y\in U^{\alpha_{0}}$,
$Y^{\alpha_{0}}_{y}=\xi$ and $y\not\in U^{\alpha}$
for $\alpha\neq\alpha_{0}$. Let $(\eta^{\alpha})$
be a partition of unity subordinated to $(U^{\alpha})$,
and define
$$ Y=\sum_{\alpha}\eta^{\alpha}Y^{\alpha}\;.$$
Then $Y_{y}=Y^{\alpha_{0}}_{y}=\xi$ and
$d\phi(Y_{z})= X_{\phi(z)}$
for any $z\in N$. Thus $Y\in L(N,\cF(\phi))$ and 
$d\phi(Y)=X$.
\end{ex}

\begin{prop} \label{prop:btcf.10}
Let $\cF$ be a locally transversely parallelizable
foliation of a connected manifold $N$,
let $s\!:N\ra M$ be a
surjective submersion with connected compact
$\cF$-saturated fibers, and let
$\bar{s}\!:N/\basF\ra M$ be the submersion
induced by $s$.
Then
the homomorphisms of Lie algebras
$$ d\baspi\!:L(N,\cF)\cap L(N,\cF(s))\lra L(N/\basF,\cF(\bar{s}))\;, $$
$$ d\bar{s}\!:L(N/\basF,\cF(\bar{s}))\lra\fX(M) $$
and
$$ ds\!:L(N,\cF)\cap L(N,\cF(s))\lra\fX(M) $$
have the lifting property.
\end{prop}

\begin{proof}
Example \ref{ex:btcf.9} implies that
$d\bar{s}\!:L(N/\basF,\cF(\bar{s}))\ra\fX(M)$ has
the lifting property. Because $ds=d\bar{s}\com d\baspi$
and $L(N,\cF)\subset L(N,\basF)$,
we have
$$ d\baspi(L(N,\cF)\cap L(N,\cF(s)))\subset L(N/\basF,\cF(\bar{s}))\;,$$
and it is sufficient to show that
$d\baspi\!:L(N,\cF)\cap L(N,\cF(s))\ra L(N/\basF,\cF(\bar{s}))$
has the lifting property.

Write $W=N/\basF$.
Take any projectable vector field
$Z\in L(W,\cF(\bar{s}))$, and choose
$\xi\in T_{y}(N)$ such that $d\baspi(\xi)=Z_{\baspi(y)}$.
Put $w=\baspi(y)$.

First we shall show that there exists an open neighbourhood
$V$ of $w$ in $W$ and a vector field
$Y^{V}\in L(\baspiinv(V),\cF|_{\baspiinv(V)})$ such that
$Y^{V}_{y}=\xi$ and $d\baspi(Y^{V})=Z|_{V}$.
To see this, choose projectable vector fields
$Y_{1},\ldots,Y_{q}\in L(N,\cF)$ such that
$y$ is a regular point of the associated
local transverse parallelism
$(\bY_{1},\ldots,\bY_{q})$ on $(N,\cF)$. 
Furthermore, we can choose these projectable
vector fields so that
$\xi=\sum_{i=1}^{q}c_{i}(Y_{i})_{y}$
for some constants $c_{1},\ldots,c_{q}$. Indeed,
we can replace $Y_{1}$ by the vector
field $Y_{1}+Y'$, for a suitable vector field
$Y'\in\fX(\cF)$, so that $\xi$ is in the span
of tangent vectors $(Y_{1})_{y},\ldots,(Y_{q})_{y}$.

Denote $X_{i}=\baspi(Y_{i})$, for $i=1,\ldots,q$.
We can reorder $Y_{1},\ldots,Y_{q}$, so that
$((X_{1})_{w},\ldots,(X_{\bar{q}})_{w})$ is a basis
of $T_{w}(W)$. Therefore we can choose an open neighbourhood
$V$ of $w$ in $W$ such that $((X_{1})_{v},\ldots,(X_{\bar{q}})_{v})$ 
is a basis of $T_{v}(W)$ for any $v\in V$.
By Example \ref{ex:btcf.9}
there exists a vector field $Z'\in L(W,\cF(\bar{s}))$
such that
$Z'_{w}=\sum_{i=1}^{\bar{q}}c_{i}(X_{i})_{w}$.
Write
$$ Z|_{V}-Z'|_{V}-\sum_{j=\bar{q}+1}^{q} c_{j} X_{j}|_{V}
   = \sum_{i=1}^{\bar{q}}h_{i} X_{i}|_{V} $$
and
$$ Z'|_{V}=\sum_{i=1}^{\bar{q}}h'_{i} X_{i}|_{V}\;, $$
so $h_{i}, h'_{i}\in \eC(V)$, $h_{i}(w)=0$
and $h'_{i}(w)=c_{i}$, for $i=1,\ldots,\bar{q}$.

Let $g_{i}=h_{i}+h'_{i}$ for $i=1,\ldots,\bar{q}$,
and put $g_{i}=c_{i}$ for $i=\bar{q}+1,\ldots,q$.
Thus for any $i=1,\ldots,q$ we have
$g_{i}\in\eC(V)$, $g_{i}(w)=c_{i}$ and
$Z|_{V}=\sum_{i=1}^{q}g_{i}X_{i}|_{V}$.
Now we define $Y^{V}\in L(\baspiinv(V),\cF|_{\baspiinv(V)})$ by
$$ Y^{V}=\sum_{i=1}^{q}(g_{i}\com\baspi)Y_{i}|_{\baspiinv(V)}\;.$$
We have
$Y^{V}_{y}=\sum_{i=1}^{q} g_{i}(w)(Y_{i})_{y}=\xi$
and
$d\baspi(Y^{V})=\sum_{i=1}^{q} g_{i} X_{i}|_{V}=Z|_{V}$.
   
From this, we can conclude that there exist an open cover
$(V^{\alpha})$ of $W$ and projectable vector fields
$Y^{\alpha}=Y^{V^{\alpha}}\in L(\baspiinv(V^{\alpha}),
\cF|_{\baspiinv(V^{\alpha})})$
such that
$d\baspi(Y^{\alpha})=Z|_{V^{\alpha}}$.
Furthermore, we can assume that $w\in V^{\alpha_{0}}$,
$Y^{\alpha_{0}}_{y}=\xi$ and $y\not\in V^{\alpha}$
for $\alpha\neq\alpha_{0}$. Let $(\eta^{\alpha})$
be a partition of unity subordinated to $(V^{\alpha})$,
and define
$$ Y=\sum_{\alpha}(\eta^{\alpha}\com\baspi) Y^{\alpha}\;.$$
Observe that $Y\in L(N,\cF)$ because the functions
$\eta^{\alpha}\com\baspi$ are basic, and that
$Y_{y}=Y^{\alpha_{0}}_{y}=\xi$.
Furthermore, 
$d\baspi(Y_{z})=\sum_{\alpha}\eta^{\alpha}(\baspi(z)) Z_{\baspi(z)})
 = Z_{\baspi(z)}$
for any $z\in N$.
Finally, we have $Y\in L(N,\cF(s))$ because
$d\baspi(Y)=Z\in L(W,\cF(\bar{s}))$.
\end{proof}

\begin{proof}[Proof of Theorem \ref{theo:btcf.11}]
Since $T_{y}(\cF)\subset T_{y}(\basF)\subset \Ker(ds)_{y}\subset T_{y}(N)$,
we can choose tangent vectors
$\xi_{1},\ldots,\xi_{q}\in T_{y}(N)$
such that their projections to $\nu_{y}(\cF)$ form a basis
of $\nu_{y}(\cF)$,
$ds(\xi_{m+1})=\ldots=ds(\xi_{q})=0$
and
$d\baspi(\xi_{\bar{q}+1})=\ldots=d\baspi(\xi_{q})=0$.
By Proposition \ref{prop:btcf.10} we can choose vector fields
$Z_{1},\ldots,Z_{\bar{q}}\in L(N/\basF,\cF(\bar{s}))$ such that
$$ (Z_{i})_{\baspi(y)}=d\baspi(\xi_{i})\;,
   \;\;\;\;\;\;\;\;\;\; i=1,\ldots,\bar{q}\;. $$
Set
$$ Z_{i}=0 \;,\;\;\;\;\;\;\;\;\;\; i=\bar{q}+1,\ldots,q\;.$$
Again by Proposition \ref{prop:btcf.10},
we can find $Y_{1},\ldots,Y_{q}\in L(N,\cF)\cap L(N,\cF(s))$
such that
$(Y_{i})_{y}=\xi_{i}$
and
$d\baspi(Y_{i})=Z_{i}$.
We can use the flows of the vector fields $Y_{1},\ldots,Y_{q}$
to obtain a local trivialization of $s$ as a bundle of
foliated manifolds, while the restrictions of
$Y_{m+1},\ldots,Y_{q}$ to a fiber $N_{s(y)}$ provide
a local transverse parallelism on $(N_{s(y)},\cF|_{N_{s(y)}})$
for which $y$ is a regular point.
\end{proof}

Theorem \ref{theo:btcf.11} motivates the following
definition. Let $M$ be a connected manifold.
A {\em bundle of transversely complete
foliations} $s\!:(N,\cF)\ra M$
with connected compact fibers over $M$
is a manifold $N$, equipped with a locally transversely
parallelizable foliation $\cF$ and a surjective
submersion $s\!:N\ra M$ with connected compact $\cF$-saturated
fibers. In particular, any such bundle
is a locally trivial fiber bundle of foliations.

For such a bundle $s\!:(N,\cF)\ra M$
we shall write
$$ L^{s}(N,\cF)=L(N,\cF)\cap \fX(\cF(s)) $$
for the Lie algebra of $s$-vertical projectable
vector fields on $(N,\cF)$,
which is the normalizer of $\fX(\cF)$
in the Lie algebra $\fX(\cF(s))$.
Let
$$ l^{s}(N,\cF)=L^{s}(N,\cF)/\fX(\cF) $$
be the associated quotient Lie algebra of
$s$-vertical transverse vector fields on $(N,\cF)$.
It is a subalgebra of the Lie algebra $l(N,\cF)$.

The subbundle $\Ker(ds)\subset T(N)$ of
$s$-vertical tangent vectors, which is the tangent bundle
of the foliation $\cF(s)$, will also be denoted by
$T^{s}(N)$.
We shall write $\nu^{s}(\cF)=T^{s}(N)/T(\cF)$
for the corresponding
subbundle of the normal bundle $\nu(\cF)$.
Note that $ds\!:T(N)\ra T(M)$ induces a map of
vector bundles $\nu(\cF)\ra T(M)$
(which we denote again by $ds$) with kernel $\nu^{s}(\cF)$.
A transversal vector field on $(N,\cF)$ is
$s$-vertical if and only if it is a section
of $\nu^{s}(\cF)$.

\section{Developability of bundles of transversely complete foliations} \label{sec:dbtcf}

Let $s\!:(N,\cF)\ra M$ be a locally trivial fiber
bundle of foliations with connected fibers.
We shall denote by $\cF_{x}$ the restriction of
$\cF$ to the fiber $N_{x}=s^{-1}(x)$ over a point $x\in M$.
We say that such a bundle of foliations $s\!:(N,\cF)\ra M$
is {\em developable}
if the foliation $(N_{x},\cF_{x})$ is developable
for any $x\in M$.

Our aim in this section is to provide some characterization
of this notion of developability, first in terms
of the fiberwise fundamental groupoid of $s$
(Theorem \ref{theo:dbtcf.5}),
and second in terms of the integrability
of an associated Lie algebroid constructed in
Subsection \ref{subsec:bla}
(Theorem \ref{theo:dbtcf.15}).

\subsection{The fiberwise fundamental groupoid} \label{subsec:ffg}

Let $s\!:(N,\cF)\ra M$ be a locally trivial fiber
bundle of foliations with connected fibers.
For any open subset $U$ of $M$, we shall write
$N|_{U}=s^{-1}(U)$ and $\cF|_{U}=\cF|_{s^{-1}(U)}$.
Note that $s|_{U}\!:(N|_{U},\cF|_{U})\ra U$ is
again a locally trivial fiber bundle of foliations.

Denote by $\Pi^{s}(N)$ the
monodromy groupoid of $(N,\cF(s))$.
This groupoid can be viewed as the fiberwise
fundamental groupoid of the bundle $s$.
The space
$\Pi^{s}(N)$ is in fact a locally trivial fiber bundle
over $M$, and the map
$(\trg,\src)\!:\Pi^{s}(N)\ra N\times_{M}N$ is a covering
projection. We define the pull-back foliation
$\Pi^{s}(\cF)=(\trg,\src)^{\ast}(\cF\times 0)$,
which is a right invariant
foliation of the groupoid $\Pi^{s}(N)$.
If $M$ is a one-point space, the Lie groupoid
$\Pi^{s}(N)$ is simply the fundamental groupoid $\Pi(N)$;
the corresponding foliation $\Pi^{s}(\cF)$
will in this case be denoted by $\Pi(\cF)$.

For a local section $\sigma\!:U\ra N$ of $s$,
defined on an open subset $U$ of $M$, we define
$\Pi^{s}_{\sigma}(N)$ by the following pull-back:
$$
\xymatrix{
\Pi^{s}_{\sigma}(N) \ar[d]_{\tau} \ar[r]
    & \Pi^{s}(N) \ar[d]^{\src} \ar[r]^-{\trg} & N\\
U \ar[r]^{\sigma} & N &
}
$$
The elements of $\Pi^{s}_{\sigma}(N)$ are
the homotopy classes of paths in $s$-fibers starting
at a point in $\sigma(U)$.
In particular, the restriction of the target map
of $\Pi^{s}(N)$ to the fiber
$\Pi^{s}_{\sigma(x)}(N)=\Pi^{s}(N)(\sigma(x),\oo)$
over a point $x\in U$
is the universal covering projection onto $N_{x}$,
and
the restriction of $\Pi^{s}(\cF)$ to 
$\Pi^{s}_{\sigma(x)}(N)$ is the pull-back
of $\cF_{x}$ along the covering projection
$\trg\!:\Pi^{s}_{\sigma(x)}(N)\ra N_{x}$.

Since $s$ is a locally trivial fiber bundle of foliations,
it follows that the map
$$ s\com\src=s\com\trg\!:(\Pi^{s}(N),\Pi^{s}(\cF))\lra M $$
is a locally trivial fiber bundle of foliated Lie groupoids, i.e.
for any $x_{0},x\in M$ there is an open neighbourhood $U$ of $x$
in $M$, and an isomorphism of Lie groupoids
$$ \Pi^{s}(N)|_{s^{-1}(U)} \lra \Pi(N_{x_{0}})\times U $$
over $U$, which maps the restriction of the foliation $\Pi^{s}(\cF)$
to the foliation $\Pi(\cF_{x_{0}})\times 0$.
In particular, the map
$\src\!:(\Pi^{s}(N),\Pi^{s}(\cF))\ra N$ is a locally
trivial fiber bundle of foliations with connected fibers.

\begin{lem} \label{lem:dbtcf.1}
Let $s\!:(N,\cF)\ra M$ be a locally trivial fiber bundle
of foliations with connected fibers.
If $\sigma\!:U\ra N$ is a local section
of $s$, defined on an open subset $U$ of $M$,
then $\trg\!:\Pi^{s}_{\sigma}(N)\ra N|_{U}$ is
a covering projection, and
$$ (\Pi^{s}_{\sigma}(N),\trg^{\ast}(\cF|_{U}))\lra U $$
is a locally trivial fiber bundle of foliated manifolds.
\end{lem}

\begin{proof}
Take any $x_{0},x_{1}\in U$, and choose an open
simply connected neighbourhood $V$
of $x_{1}$ in $U$ such that $(N|_{V},\cF|_{V})$ is a 
trivial bundle of foliated manifolds, with
a trivialization
$$ \mu\!=(\delta,s):N|_{V}\lra N_{x_{0}}\times V $$
which maps the $\cF|_{V}$ to the foliation
$\cF_{x_{0}}\times 0$. Choose $y_{0}\in N_{x_{0}}$.
Since $\trg\!:\Pi^{s}_{y_{0}}(N)\ra N_{x_{0}}$ is a covering
projection and $V$ is simply connected, we can lift
$\delta\com\sigma|_{V}$ to a map
$\tau\!:V\ra \Pi^{s}_{y_{0}}(N)$, so
$\tau(x)\!:y_{0}\ra \delta(\sigma(x))$ for any $x\in V$.

Define a map
$\tilde{\mu}\!:\Pi^{s}_{\sigma|_{V}}(N)\ra
  \Pi^{s}_{y_{0}}(N)\times V$
by
$$ \tilde{\mu}(\alpha)=( \delta_{\ast}(\alpha)\tau(x),x) \;,$$
for any $\alpha\!:\sigma(x)\ra y$ in $\Pi^{s}_{\sigma|_{V}}(N)$.
This map is clearly a diffeomorphism, and the diagram
$$
\xymatrix{
\Pi^{s}_{\sigma|_{V}}(N) \ar[r]^-{\tilde{\mu}} \ar[d]^{\trg} &
    \Pi^{s}_{y_{0}}(N)\times V
    \ar[d]^{\trg\times\id} \\
N|_{V} \ar[r]^-{\mu} \ar[d]^{s} & N_{x_{0}}\times V \ar[dl]^{\pr_{2}} \\
V & 
}
$$
commutes. This shows that 
$\trg\!:\Pi^{s}_{\sigma}(N)\ra N|_{U}$ is
a covering projection. Since
the map $\mu$ trivializes the bundle
$(N|_{V},\cF|_{V})\ra V$ of foliated manifolds, it follows that
$(\Pi^{s}_{\sigma}(N),\trg^{\ast}(\cF|_{U}))\ra U$
is also a locally trivial fiber bundle of foliated manifolds.
\end{proof}

We first observe that strict simplicity is a fiberwise property:

\begin{prop} \label{prop:dbtcf.3}
Let $s\!:(N,\cF)\ra M$ be a locally trivial fiber bundle
of foliations with connected fibers. Then
the foliation $\cF$ is strictly simple if and only if
for each $x\in M$ the foliation
$\cF_{x}$ is strictly simple.
\end{prop}

\begin{proof}
Suppose that $\cF$ is strictly simple, thus given by
a submersion $f\!:N\ra T$ with connected fibers.
Since the fibers of $s$ are saturated, the map
$f$ induces a surjective
submersion $\bar{f}\!:T\ra M$. In particular,
the fiber $T_{x}=\bar{f}^{-1}(x)$ is a closed submanifold of $T$.
Moreover, the restriction
$$ f|_{N_{x}}\!:N_{x}\lra T_{x} $$
is a submersion with connected fibers and defines the
foliation $\cF_{x}$.

Conversely, suppose that
$\cF_{x}$ is strictly simple for any $x\in M$.
Equivalently, there is a structure of a Hausdorff manifold
on the space of leaves $N_{x}/\cF_{x}$
such that the quotient projection
$N_{x}\ra N_{x}/\cF_{x}$ is a submersion.
Since $N$ is a locally trivial fiber bundle
of foliated manifolds, 
$N/\cF$ can also be given a
structure of a Hausdorff manifold such that the
quotient projection
$N\ra N/\cF$
is a submersion. It follows that $\cF$ is strictly
simple.
\end{proof}

\begin{theo} \label{theo:dbtcf.5}
Let $s\!:(N,\cF)\ra M$ be a locally trivial fiber bundle
of foliations with connected fibers. The following
conditions are equivalent.

(i)
The bundle $s\!:(N,\cF)\ra M$ is developable.

(ii)
The foliation $\Pi^{s}(\cF)$ of $\Pi^{s}(N)$
is strictly simple.

(iii)
For any open subset $U$ of $M$ with a section $\sigma\!:U\ra M$ of $s$,
the restriction of the foliation $\Pi^{s}(\cF)$ to $\Pi^{s}_{\sigma}(N)$
is strictly simple.
\end{theo}

\begin{proof}
(i)$\Leftrightarrow$(ii)
Recall that
$\src\!:(\Pi^{s}(N),\Pi^{s}(\cF))\ra N$ is a locally
trivial fiber bundle of foliations with connected fibers.
If we apply Proposition \ref{prop:dbtcf.3}
to this bundle, we see that
$\Pi^{s}(\cF)$ is strictly simple if and only if
the restriction of $\Pi^{s}(\cF)$ to
$\Pi^{s}(N)(y,\oo)$ is strictly simple, for any $y\in N$.
But
$$ \trg\!:\Pi^{s}(N)(y,\oo)\lra N_{s(y)} $$
is the universal
cover of $N_{s(y)}$, and the restriction of
$\Pi^{s}(\cF)$ to $\Pi^{s}(N)(y,\oo)$ is the pull-back of
$\cF_{s(y)}$ along this covering.

(i)$\Leftrightarrow$(iii)
By definition, the bundle
$s\!:(N,\cF)\ra M$ is developable
if and only if the pull-back of the foliation $(N_{x},\cF_{x})$
to the universal cover of $N_{x}$ is strictly simple,
for any $x\in M$. This is true if and only if the foliation
$(\Pi^{s}_{\sigma(x)}(N),\trg^{\ast}(\cF_{x}))$ is strictly
simple, for any section $\sigma\!:U\ra N$ of $s$, because
$\Pi^{s}_{\sigma(x)}(N)$ is the universal cover of $N_{x}$.
In turn, this holds true
if and only if the foliation $(\Pi^{s}_{\sigma}(N),\trg^{\ast}(\cF|_{U}))$
is strictly simple. Indeed, to see this, apply
Proposition \ref{prop:dbtcf.3} to the map
$$ s\com\trg\!:(\Pi^{s}_{\sigma}(N),\trg^{\ast}(\cF|_{U}))\lra U\;,$$
which is a locally trivial fiber bundle of foliated manifolds.
Since $\Pi^{s}_{\sigma}(N)$ is a covering space of $N|_{U}$ by
Lemma \ref{lem:dbtcf.1}, this completes the proof.
\end{proof}

\begin{cor} \label{cor:dbtcf.55}
For a foliation $\cF$ on a manifold $M$, the following two
statements are equivalent:
\begin{enumerate}
\item [(i)]  The foliation $\cF$ is developable.
\item [(ii)] The foliation $\cF$, viewed as a subalgebroid of $T(M)$,
             is integrable by a closed subgroupoid of the fundamental
             groupoid $\Pi(M)$.
\end{enumerate}
\end{cor}

\begin{proof}
The groupoid $\Pi(M)$ integrates the algebroid $T(M)$, so
the subalgebroid $\cF\subset T(M)$ correspond to 
an invariant foliation of $\Pi(M)$ by
\cite[Lemma 2.1]{MoerdijkMrcun2004a},
which is readily identified to the foliation $\Pi(\cF)$.
Thus by Theorem \ref{theo:dbtcf.5},
$\cF$ is developable if and only if $\Pi(\cF)$ is
strictly simple, and by
\cite[Proposition 3.2 (ii)]{MoerdijkMrcun2004a}
this is the case if and only if $\cF$ can be integrated
by a closed subgroupoid of $\Pi(M)$.
\end{proof}

\subsection{The basic Lie algebroid} \label{subsec:bla}

As the main result of this section,
we shall now construct
a Lie algebroid associated
to a bundle $s\!:(N,\cF)\ra M$ of transversely
complete foliations, and
show that integrability of this algebroid
is equivalent to developability of the bundle $s$.

\begin{lem} \label{lem:dbtcf.6}
Let $s\!:(N,\cF)\ra M$ be a bundle of transversely
complete foliations with compact connected fibers,
and let $\baspi\!:N\ra N/\basF=W$ be the associated
basic fibration.
The groupoid $N\times_{W}N$ naturally acts 
linearly both on
the normal bundle $\nu(\cF)$ and its subbundle
$\nu^{s}(\cF)$ along $\baspi$.
\end{lem}

\begin{proof}
For any
$(z,y)\in N\times_{W}N$ and $\xi\in\nu_{z}(\cF)$,
put
$$ \Theta(\xi,z,y)=\bY_{y} \;,$$
where $\bY$ is any transverse vector field on $(N,\cF)$
such that $\bY_{z}=\xi$.
This definition is independent of the
choice of $\bY$ by Remark \ref{rems:btcf.1} (5).
Using local transverse parallelisms and the corresponding local
trivializations of the normal bundle,
it is straightforward to check
that $\Theta$ is indeed a smooth linear action.

This action restricts to an action of $N\times_{W}N$
on $\nu^{s}(\cF)$ by
Proposition \ref{prop:btcf.10}.
Indeed,
for any $(z,y)\in N\times_{W}N$ and any $\xi\in\nu^{s}_{z}(\cF)$,
we can find $\bY\in L^{s}(N,\cF)$ with
$\bY_{z}=\xi$.
\end{proof}

\begin{rem} \rm \label{rem:dbtcf.7}
Note that if $y$ and $z$ lie on the same leaf of $\cF$,
the corresponding linear map
$\Theta(\oo,z,y)\!:\nu_{z}(\cF)\ra\nu_{y}(\cF)$
is simply the linear holonomy
isomorphism of an arbitrary path in the leaf of
$\cF$ from $z$ to $y$.
\end{rem}

Let $s\!:(N,\cF)\ra M$ be a bundle of transversely
complete foliations with compact connected fibers
over a connected manifold $M$, and write $W=N/\basF$.
We shall denote by
$$ \fb^{s}(N,\cF)=\nu^{s}(\cF)/(N\times_{W}N) $$
the space of orbits of the natural action of
the groupoid $N\times_{W}N$ on $\nu^{s}(\cF)$
(Lemma \ref{lem:dbtcf.6}).
In the next lemma, we
show that the natural projection
$\fb^{s}(N,\cF)\ra W$ is a vector bundle.

\begin{lem} \label{lem:dbtcf.8}
Let $s\!:(N,\cF)\ra M$ be a bundle of transversely
complete foliations with compact connected fibers, and let
$\baspi\!:N\ra N/\basF=W$ be the associated basic fibration.
Then the orbit space $\fb^{s}(N,\cF)$
has the structure of a vector bundle over $W$
such that

(i)
the quotient projection
$\theta\!:\nu^{s}(\cF)\ra\fb^{s}(N,\cF)$
is a $(N\times_{W}N)$-principal bundle, and

(ii)
the diagram
$$
\xymatrix{
\nu^{s}(\cF) \ar[d] \ar[r]^-{\theta} & \fb^{s}(N,\cF) \ar[d] \\
N \ar[r]^{\baspi} & W
}
$$
is a fibered product of vector bundles.
\end{lem}

\begin{proof}
Since $N$ is a principal $N\times_{W}N$-bundle over $W$
and the map $\nu^{s}(\cF)\ra N$ is $N\times_{W}N$-equivariant,
it follows that there is a structure of smooth Hausdorff
manifold on
$\fb^{s}(N,\cF)$
such that the quotient projection
$\theta$
is a $N\times_{W}N$-principal bundle.
Since the action on $\nu^{s}(\cF)$ is given
by linear isomorphisms, the fibers of $\fb^{s}(N,\cF)\ra W$
have vector
space structures such that the restriction
$\theta_{y}\!:\nu^{s}_{y}(\cF)\ra\fb^{s}(N,\cF)_{\baspi(y)}$
of $\theta$ is a linear isomorphism, for any $y\in N$.
Indeed, if $z\in N$ is another point on the same basic leaf as $y$,
then $\theta_{y}\com\Theta(\oo,z,y)=\theta_{z}$,
where $\Theta$ denotes the natural action of
$N\times_{W}N$ on $\nu(\cF)$.

Take any local transverse parallelism
$(\bY_{1},\ldots,\bY_{q})$ with regular set $U$,
which is of the form given by Theorem \ref{theo:btcf.11}.
It gives us a local trivialization of $\fb^{s}(N,\cF)$
$$ \alpha\!:\baspi(U)\times\RR^{q-m}\lra
   \fb^{s}(N,\cF)|_{\baspi(U)} $$ 
over $\baspi(U)$ by
$$ \alpha(\baspi(y),t_{m+1},\ldots,t_{q})=
   \theta_{y}\left(\sum_{j=m+1}^{q}t_{j}(\bY_{j})_{y}\right)\;.$$
This is well-defined because for another point
$z\in U$ with
$\baspi(y)=\baspi(z)$ we have
$\Theta((\bY_{j})_{z},z,y)=(\bY_{j})_{y}$.
Thus we can conclude that $\fb^{s}(N,\cF)$ is indeed a vector bundle over
$W$. The diagram in (ii) is a pull-back
because $\theta$ restricts to the isomorphism $\theta_{y}$
on each fiber $\nu^{s}_{y}(\cF)$.
\end{proof}

\begin{prop} \label{prop:dbtcf.12}
Let $s\!:(N,\cF)\ra M$ be a bundle of transversely
complete foliations with compact connected fibers.
The bundle $\fb^{s}(N,\cF)$
has a natural
structure of a regular Lie algebroid over $N/\basF$
such that

(i)
the quotient projection
$\nu^{s}(\cF)\ra \fb^{s}(N,\cF)$
over the basic fibration induces
an isomorphism of Lie algebras
$$ \Gamma(\fb^{s}(N,\cF))\lra l^{s}(N,\cF)\;,$$

(ii)
the foliation of $N/\basF$
corresponding to this regular algebroid
is the simple foliation given by
the submersion $N/\basF\ra M$ induced by $s$, and

(iii)
the natural projection
$T^{s}(N) \ra \fb^{s}(N,\cF)$
is a surjective morphism of Lie algebroids.
\end{prop}

\begin{proof}
The fact that the diagram
in Lemma \ref{lem:dbtcf.8} (ii)
is a fibered product
implies that any section $X$ of $\fb^{s}(N,\cF)$
induces a section
$\baspiast X$ of $\nu^{s}(\cF)$, so we get a map
$$ \baspiast\!:\Gamma(\fb^{s}(N,\cF))\lra\Gamma(\nu^{s}(\cF))\;,$$
which is clearly $\basOmega(N,\cF)$-linear
(recall that the composition with the
basic projection $\baspi\!:N\ra N/\basF=W$
gives us the
identification $\basOmega(N,\cF)=\eC(W)$).

Note that a section $\sigma$ of the bundle
$\nu^{s}(\cF)$ is a transverse vector field if and only if
it is holonomy invariant, i.e.\ if and only if
\begin{equation}  \label{eq:dbtcf.1}
\Theta(\sigma_{z},z,y)=\sigma_{y}
\end{equation}
for any two points $y,z$ on the same leaf of $\cF$.
Furthermore, 
Remark \ref{rems:btcf.1} (5) implies
that for a transverse
vector field the condition (\ref{eq:dbtcf.1})
in fact automatically
holds for any $y,z$ on the same leaf of $\basF$.
Thus, by the definition of $\fb^{s}(N,\cF)$ and
by Lemma \ref{lem:dbtcf.8},
the $s$-vertical transverse vector fields
on $(N,\cF)$
are exactly those sections
of $\nu^{s}(\cF)$ which can be projected along
$\baspi$ to a section
of $\fb^{s}(N,\cF)$. We can therefore conclude that
$$ \baspiast\!:\Gamma(\fb^{s}(N,\cF))\lra l^{s}(N,\cF) $$
is an $\basOmega(N,\cF)$-linear isomorphism.
We therefore define a Lie algebra structure on
$\Gamma(\fb^{s}(N,\cF))$
so that $\baspiast$ is also a Lie algebra isomorphism.

To give $\fb^{s}(N,\cF)$ a Lie algebroid structure, we have to define
its anchor, which is a morphism of vector bundles
$\anchor\!:\fb^{s}(N,\cF)\lra T(W)$
over $W$. We define it by
$$ \anchor(\theta_{y}(\zeta+T_{y}(\cF)))=(d\baspi)_{y}(\zeta) $$
for any $\zeta\in T^{s}_{y}(N)$.
It is straightforward to check that
$\fb^{s}(N,\cF)$ is a
Lie algebroid over $W$
and that $T^{s}(N) \ra \fb^{s}(N,\cF)$ is
a surjective morphism of Lie algebroids.
\end{proof}

The Lie algebroid
$$ \fb^{s}(N,\cF) $$
over $N/\basF$
will be called the
{\em basic Lie algebroid}
associated to the bundle of transversally complete foliations
$s\!:(N,\cF)\ra M$.

\begin{theo} \label{theo:dbtcf.15}
Suppose that $s\!:(N,\cF)\ra M$ is a bundle
of transversely complete
foliations with compact connected fibers.
Then $s\!:(N,\cF)\ra M$ is developable if and only if
the associated basic Lie
algebroid $\fb^{s}(N,\cF)$ is integrable.
\end{theo}

\begin{proof}
Recall from Theorem \ref{theo:btcf.11} that the restriction
$\cF_{x}$ of $\cF$ to the fiber $N_{x}$ over $x\in M$
is a transversely complete
(and therefore homogeneous) foliation.
Theorem  \ref{theo:btcf.11}
also implies $(\cF_{x})_{\mathrm{bas}}=\basF|_{N_{x}}$,
so the fiber $W_{x}=\bar{s}^{-1}(x)$ over $x$
is the space of leaves $N_{x}/(\cF_{x})_{\mathrm{bas}}$.

Since $\fb^{s}(N,\cF)$ is a locally trivial bundle
of transitive Lie algebroids over $M$, its integrability is
equivalent to the integrability of each fiber.
The restriction of $\fb^{s}(N,\cF)$ to the fiber $W_{x}$ is
exactly the basic algebroid $\fb(N_{x},\cF_{x})$
of the transversely complete foliation $\cF_{x}$
of the compact fiber $N_{x}$ (see also \cite{Molino1988}).
On the other hand, developability of $s\!:(N,\cF)\ra M$
is also given fiberwise. Therefore we may
assume without loss of generality that $M$ is a one-point space.
In this case, the statement follows from the
Almeida-Molino theorem for transversely complete
foliations of compact manifolds \cite{Molino1988}
(see also \cite{MoerdijkMrcun2003}).
\end{proof}

\section{Proof of the main theorem} \label{sec:pmt}

In this final
section, we return to the setting of Section \ref{sec:mcfds},
and use the results of Sections
\ref{sec:btcf} and \ref{sec:dbtcf} to prove
Theorem \ref{theo:mcfds.25}.
In the proof of this theorem,
we use the following simple observations on the
universal covering groupoid $\tilde{G}$ and the fiberwise
fundamental groupoid $\Pi^{\src}(G)$.

\begin{lem} \label{lem:mcfds.30}
Let $G$ be a source-connected Lie groupoid
with Lie algebroid $\fg$.

(i)
There is a natural action of $G$ on $\Pi^{\src}(G)$
which
differentiates to an action $\nabla$ by $\fg$
on the Lie algebroid $T^{\src}(G)\ra G$ along the source map
$\src\!:G\ra M$. For any $X\in\Gamma(\fg)$, the corresponding
derivation $\nabla(X)=(\nabla_{X},R(X))$ is given by
$$ R(X)=\tilde{X}^{-1} $$
and
$$ \nabla_{X}(Y)=[R(X),Y] $$
for any $Y\in\fX^{\src}(G)$,
where $\tilde{X}^{-1}$ the image of
the right invariant extension $\tilde{X}$
of $X$ along the inverse
map of $G$.

(ii)
The universal covering groupoid $\tilde{G}$ of $G$
embeds into $G \ltimes \Pi^{\src}(G)$
as the full subgroupoid on the space of units
$\uni\!:M\ra G$, and this inclusion differentiates
to the diagonal embedding $\fg\ra\fg\ltimes T^{\src}(G)$
of Lie algebroids.
\end{lem}

\begin{proof}
Note that $\Pi^{\src}(G)$ is a bundle of groupoids over $M$
whose fiber over $x\in M$ is the fundamental groupoid
$\Pi(\src^{-1}(x))$ of
$\src^{-1}(x)$. The natural left action
of an arrow $g\!:x\ra y$ of $G$ is the map
of fundamental groupoids $\Pi(\src^{-1}(x))\ra\Pi(\src^{-1}(y))$
induced by the right translation $R_{g^{-1}}$.
The semi-direct product $G\ltimes\Pi^{\src}(G)$
can be described explicitly as the groupoid over $G$
whose arrows $g\ra g'$ are pairs $(h,\alpha)$,
where $h\!:\src(g)\ra\src(g')$ is an arrow in $G$ and
$\alpha\!:g\ra g'h$ is an arrow in $\Pi^{\src}(G)$.
From this description it is clear that the restriction
of $G\ltimes\Pi^{\src}(G)$ to the units is precisely $\tilde{G}$.
The rest follows by straightforward inspection.
\end{proof}

\begin{proof}[Proof of Theorem \ref{theo:mcfds.25}]
Before proving (i)-(iv), let us observe that
since the source map of $G$ has compact fibers
and the foliation $\cF(\fh)\subset\cF(\src)$
is transversely complete, the map
$\src\!:(G,\cF(\fh))\ra M$ is a bundle of
transversely complete foliations with connected compact
fibers. Moreover, the foliations
$\cF(\fh)$ and $\cF(\src)$, as well as the
basic foliation
$\cF(\fh)_{\mathrm{bas}}$ associated to $\cF(\fh)$,
are all right invariant. 

The closure $\bar{H}$ of $H$ in $G$ is a closed
subgroupoid of $G$ which corresponds to the basic foliation
$\cF(\fh)_{\mathrm{bas}}$, and
the space of cosets $G/\bar{H}$
is the associated space of basic leaves
\cite[Theorem 3.7]{MoerdijkMrcun2004a}.
Furthermore, the quotient
projection $\pi=\baspi\!:G\ra G/\bar{H}$ is
a locally trivial fiber bundle.

(i)
We take $\fb(G,\fh)$ to be the basic Lie algebroid
$\fb^{\src}(G,\cF(\fh))$ associated to the
bundle $\src\!:(G,\cF(\fh))\ra M$.
Proposition \ref{prop:dbtcf.12} (iii) implies
that we have
a natural surjective morphism of Lie algebroids
$\omega\!:T^{\src}(G) \ra \fb(G,\fh)$
with kernel $\cF(\fh)$.
\begin{equation}  \label{eq:mcfds.41}
\xymatrix{
T^{\src}(G) \ar[d] \ar[r]^-{\omega} & \fb(G,\fh) \ar[d] \\
G \ar[r]^-{\pi} & G/\bar{H}
}
\end{equation}

(ii)
Note that the construction
of $\fb(G,\fh)$ in Section \ref{sec:dbtcf}
is invariant under the right $G$-action,
so $G$ acts on the Lie algebroids $T^{\src}(G)$
and $\fb(G,\fh)$, and the map $\omega$
is invariant under this action.
This action, formally inverted to a left action,
differentiates to an action of $\fg$ on $\fb(G,\fh)$.

(iii)
We claim that the restriction of  $\omega$
to $\fg\subset T^{\src}(G)$ provides
a Maurer-Cartan form with kernel $\cF(\fh)$.
By the action of $\fg$ on $\fb(G,\fh)$,
diagram (\ref{eq:mcfds.41}) gives us a map of Lie algebroids
over $G$,
$$ \fg\ltimes\omega\!: \fg\ltimes T^{\src}(G) \lra
   \fg\ltimes\fb(G,\fh) \;,$$
with kernel $\fg\ltimes\cF(\fh)$.
Precomposing this map with the natural
diagonal section $\fg\ra\fg\ltimes T^{\src}(G)$ over
the unit section $\uni\!:M\ra G$,
Lemma \ref{lem:mcfds.30} provides
a morphism of Lie algebroids
$$ \fg\lra\fg\ltimes\fb(G,\fh) $$
over $\pi\com\uni\!:M\ra G/\bar{H}$,
which gives us the required Maurer-Cartan form.

(iv)
The foliation $\cF(\fh)$ of $G$ by cosets of $H$
is locally transversely parallelizable, and
by Theorem \ref{theo:btcf.11} 
the source map $\src\!:(G,\cF(\fh))\ra M$ is
a locally trivial bundle of transversely complete
foliations with compact connected fibers.
By Theorem \ref{theo:dbtcf.15} the
basic Lie algebroid $\fb(G,\fh)$
is integrable if and only if
the bundle $\src\!:(G,\cF(\fh))\ra M$ is developable,
and by Theorem \ref{theo:dbtcf.5}
this is true if and only  if
the corresponding foliation 
$\Pi^{\src}_{\uni}(\cF(\fh))$ of $\tilde{G}=\Pi^{\src}_{\uni}(G)$
is strictly simple
(here $\uni$ is the unit section of the source map).
By Proposition \ref{prop:mcfds.4} this is equivalent
to developability of the Lie algebroid $\fh$.
\end{proof}

\begin{rem} \rm \label{rem:mcfds.27}
The vector bundle $\fb(G,\fh)$ is a
locally trivial bundle of Lie algebroids over $M$,
the fiber over $x\in M$ being a transitive
algebroid over $\src^{-1}(x)\subset G/\bar{H}$.
If $\fb(G,\fh)$ is integrable,
write $B$ for its source-simply connected integral.
Then $B$ is a locally trivial bundle of transitive
Lie groupoids over $M$. It follows that
$B$ is Hausdorff. The action of $G$ on
$\fb(G,\fh)$ integrates to an action
of $G$ on $B$ (see \cite{MoerdijkMrcun2002}),
and the Maurer-Cartan form integrates to a twisted homomorphism
$F\!:\tilde{G}\ra B$, whose kernel is the closed subgroupoid
integrating $\fh$.
\end{rem}


\begin{thebibliography}{10}

\bibitem{AlmeidaMolino1985}
R.~Almeida and P.~Molino.
\newblock Suites d'{A}tiyah et feuilletages transversalement complets.
\newblock {\em C. R. Acad. Sci. Paris S\'er. I Math.}, 300(1):13--15, 1985.

\bibitem{BursztynWeinstein2004}
H.~Bursztyn and A.~Weinstein.
\newblock Poisson geometry and {M}orita equivalence.
\newblock {\em London Math. Soc. Lecture Note Ser.}, to appear, preprint
  arXiv:math.SG/0402347, 2004.

\bibitem{CannasdasilvaWeinstein1999}
A.~Cannas~da Silva and A.~Weinstein.
\newblock {\em Geometric models for noncommutative algebras}, volume~10 of {\em
  Berkeley Mathematics Lecture Notes}.
\newblock American Mathematical Society, Providence, RI, 1999.

\bibitem{CattaneoFelder2000}
A.~S. Cattaneo and G.~Felder.
\newblock A path integral approach to the {K}ontsevich quantization formula.
\newblock {\em Comm. Math. Phys.}, 212(3):591--611, 2000.

\bibitem{Conlon1974}
L.~Conlon.
\newblock Transversally parallelizable foliations of codimension two.
\newblock {\em Trans. Amer. Math. Soc.}, 194:79--102, 1974.

\bibitem{Connes1994}
A.~Connes.
\newblock {\em Noncommutative geometry}.
\newblock Academic Press Inc., San Diego, CA, 1994.

\bibitem{DouadyLazard1966}
A.~Douady and M.~Lazard.
\newblock Espaces fibr\'es en alg\`ebres de {L}ie et en groupes.
\newblock {\em Invent. Math.}, 1:133--151, 1966.

\bibitem{Fedida1971}
E.~Fedida.
\newblock Sur les feuilletages de {L}ie.
\newblock {\em C. R. Acad. Sci. Paris S\'er. A-B}, 272:A999--A1001, 1971.

\bibitem{HigginsMackenzie1990}
P.~J. Higgins and K.~Mackenzie.
\newblock Algebraic constructions in the category of {L}ie algebroids.
\newblock {\em J. Algebra}, 129(1):194--230, 1990.

\bibitem{Mackenzie1987}
K.~Mackenzie.
\newblock {\em Lie groupoids and {L}ie algebroids in differential geometry},
  volume 124 of {\em London Mathematical Society Lecture Note Series}.
\newblock Cambridge University Press, Cambridge, 1987.

\bibitem{MoerdijkMrcun2002}
I.~Moerdijk and J.~Mr{\v{c}}un.
\newblock On integrability of infinitesimal actions.
\newblock {\em Amer. J. Math.}, 124(3):567--593, 2002.

\bibitem{MoerdijkMrcun2003}
I.~Moerdijk and J.~Mr{\v{c}}un.
\newblock {\em Introduction to foliations and {L}ie groupoids}, volume~91 of
  {\em Cambridge Studies in Advanced Mathematics}.
\newblock Cambridge University Press, Cambridge, 2003.

\bibitem{MoerdijkMrcun2004a}
I.~Moerdijk and J.~Mr{\v{c}}un.
\newblock On the integrability of subalgebroids.
\newblock {\em Preprint}, arXiv:math.DG /0406558, 2004.

\bibitem{Molino1977}
P.~Molino.
\newblock \'{E}tude des feuilletages transversalement complets et applications.
\newblock {\em Ann. Sci. \'Ecole Norm. Sup. (4)}, 10(3):289--307, 1977.

\bibitem{Molino1988}
P.~Molino.
\newblock {\em Riemannian foliations}, volume~73 of {\em Progress in
  Mathematics}.
\newblock Birkh\"auser Boston Inc., Boston, MA, 1988.
\newblock Translated from the French by Grant Cairns, With appendices by
  Cairns, Y. Carri\`ere, \'E. Ghys, E. Salem and V. Sergiescu.

\bibitem{Pradines1967}
J.~Pradines.
\newblock Th\'eorie de {L}ie pour les groupo\"\i des diff\'erentiables.
  {C}alcul diff\'erenetiel dans la cat\'egorie des groupo\"\i des
  infinit\'esimaux.
\newblock {\em C. R. Acad. Sci. Paris S\'er. A-B}, 264:A245--A248, 1967.

\end{thebibliography}
\end{document}